\documentclass[12pt]{article}

\usepackage{amsfonts, amssymb, amsmath, amsthm, epsf, graphicx}

\usepackage{fullpage}

\newcommand{\cP}{{\mathcal P}}
\newcommand{\cR}{{\mathcal R}}

\newcommand{\cW}{{\mathcal W}}
\newcommand{\ux}{{\underline  x}}
\newcommand{\ox}{{\overline x}}

\newcommand{\ou}{{\overline u}}
\newcommand{\uu}{{\underline u}}
\newcommand{\uU}{{\underline U}}

\newcommand{\bbN}{{\mathbb N}}
\newcommand{\bbR}{{\mathbb R}}

\newtheorem{theorem}{Theorem}[section]
\newtheorem{lemma}{Lemma}[section]

\theoremstyle{definition}
\newtheorem{definition}{Definition}[section]
\newtheorem{remark}{Remark}[section]

\numberwithin{equation}{section} \numberwithin{figure}{section}

\title{Pattern formation in parabolic equations containing hysteresis with diffusive thresholds}

\author{Pavel Gurevich\footnote{Free University Berlin, Germany \& Peoples' Friendship University of Russia, Russia; email: gurevichp@gmail.com},\qquad
Dmitrii Rachinskii\footnote{Department of Mathematical Sciences,
University of Texas at Dallas, USA \& Department of Applied
Mathematics, University College Cork, Ireland; email:
dmitry.rachinskiy@utdallas.edu}}

\begin{document}

\maketitle


\begin{abstract}
We consider a reaction-diffusion system with discontinuous reaction terms modeled by non-ideal relays.
The system is motivated by an epigenetic population model of evolution of two-phenotype bacteria
which switch phenotype in response to variations of environment. We prove the formation of patterns in the
phenotype space. The mechanism responsible for pattern formation is based
on memory (hysteresis) of the non-ideal relays.
\end{abstract}

%
%


\section{Introduction}
\begin{figure}[ht]
{\ \hfill\epsfxsize60mm\epsfbox{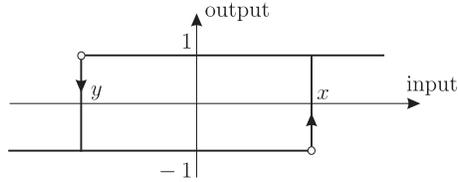}\hfill\ }
\caption{Bi-stable switch (non-ideal relay).} \label{hyst}
\end{figure}

Bi-stability is common in living systems. Max Delbr\"uck was the
first to associate different stable stationary states with
different phenotypes, or epigenetic differences, in clonal
populations such as those arising in the process of cell
differentiation \cite{delbr}. Starting from early experimental
work of Novick and Wiener on {\em lac-operon} in {\em E.~coli}, it
has been believed that bi-stability can help organisms gain
fitness in varying environmental conditions if the organism can
switch phenotype following the environmental change. {\em
Lac-operon} is a collection of genes associated with transport and
metabolism of lactose in the bacterium. Expression of these genes
can be turned on by molecules that have been called inducers.
Novick and Weiner \cite{f0} as well as Cohn and Horibata
\cite{cohn1,cohn2,cohn3}, relying on prior work of others
\cite{monod,benzer,Spiegelman}, demonstrated that two phenotypes
each associated with ``on'' and ``off'' state of {\em lac-operon}
expression can be obtained from the same culture of genetically
identical bacteria. The fraction of each corresponding
sub-population depended on the history of exposure to the inducer.
As illustrated in Fig.~\ref{hyst},
the {\em lac-operon} state was induced (switched on) when the
extracellular inducer concentration (input) exceeded an upper
threshold. The operon was switched off when the inducer
concentration fell below the lower threshold. Both phenotypes
remain stable through multiple generations of the bacterial
culture after the extracellular concentration of the inducer is
reduced to lower levels. Novick and Weiner did not use the term
``bi-stable switch'' (non-ideal relay \cite{SH}) to describe their
observations, but effectively that is what it was. Difference in
the reproductive rates of different phenotypes in different
environments was also noted in their work suggesting association
of phenotype bi-stability and switching with strategies for
increasing population growth in varying environments.

These early findings on the bi-stability of the {\em lac-operon}
 were consistent with earlier findings on stability of
 enzymes in yeast \cite{winge}.  Recent experiments using
molecular biology methods (such as those incorporating green
fluorescent protein expression under the {\em lac-operon}
promoter) permitted to confirm and further study the region of
bi-stability of the {\em lac-operon} even when multiple input
variables (TMG that acts as the inducer and glucose, for example)
were used \cite{oud}. Furthermore, multi-stable gene expression
and phenotype switching in varying environmental conditions have
been well documented in many natural and artificially constructed
biological systems
\cite{wanga,lai,grazi,f3,bac4,smits,maamar,ark,c10,15,10}.
These experiments were complemented by a substantial body of 
theoretical work that  demonstrated
that organisms can gain fitness by adjusting the rate of switching
between phenotypes, or switching thresholds, to the law of
environmental variations. A typical setting for the analytical
work has been an optimization problem where fitness (measured, for
example, by the average growth rate of the population) is the
utility function,  each phenotype is optimal for a certain state
of environment, and there is a cost associated with changing
phenotype. This setting has been applied to models based on
differential equations of population growth, switching systems and
game theory \cite{bac2, c12, h2}. One well documented experimental
fact, which can translate into the cost of phenotype switching, is
the so-called lag phase that delays reproduction in an organism
undergoing a switch. For example, in a simple (stochastic)
differential equation model incorporating the lag phase, which was
proposed in \cite{1}, the growth rate is maximized by a
bi-stability range of a certain size that depends on the
parameters of the environmental input. Another related idea that
has been explored in a similar modeling context is bet-hedging,
that is a strategy promoting diversification whereby co-existence
of several phenotypes, including those less fit to the current
environmental conditions, can favor growth in varying environments
\cite{bac1}.

The bi-stable switch (also known as the non-ideal relay) is the
simplest model of a bi-stable system. In particular, the
differential model of phenotype switching proposed in \cite{1}
reduces to the non-ideal relay in the limit of fast switching. In
the present work, we are interested in the evolutionary process in
a colony of two-phenotype bacteria where the competition for a
resource (nutrients in this case) acts as a selection pressure.
Motivated by the fact that different bi-stability ranges are
optimal (ensure the fastest growth rate) for different laws of
variation of nutrients concentration, we model bacteria by
non-ideal relays, which have all possible bi-stability ranges from
a certain interval. Moreover, we allow a diffusion process
over the interval of admissible bi-stability ranges.
We assume that each bacterium can
sporadically change its switching threshold
in response to noise in the environment or some internal processes; also, a bacterium produces offsprings with different
thresholds. Further, these random changes of the threshold
are assumed to occur according to the Gaussian
distribution. Hence, diffusion
imitates random variations of the bi-stability range that
diversify switching thresholds  in bacteria.
In sum, the evolution is
defined by two factors: competition and diffusive (random) changes
of the bi-stability range. The resulting model describes dynamics
of populations by a reaction-diffusion system with non-ideal
relays continuously distributed over an interval of bi-stability
ranges; the system also includes integral terms in the equations
for nutrients, describing the influence of bacteria on the
environment.

The well-posedness of this model was proved in \cite{2}. It was
also shown that each solution converges to a steady state
distribution (due to the diffusion, the concentration of bacteria
distributes uniformly over all the interval of admissible bi-stability
ranges). The model naturally admits a continual set of steady
state solutions as there is a continuum of bi-stable elements. Moreover, a numerical evidence was given that
in case of slow diffusion solutions converge to steady states that
can be naturally interpreted as patterns in the phenotype space.

Differential equations with non-ideal relays have been previously
used to model reaction and diffusion processes in spatially
distributed colonies of bacteria
\cite{Jaeger1,Kopfova,VisintinSpatHyst,GurevichTikhomirovEquadiff13}.
In contrast to our setting, all relays in these models have been
assumed identical (that is, having the same bi-stability range),
while spatial patterns appeared, in particular, due to the spatial
diffusion of environmental substances.

The mechanism of pattern formation in the phenotype space
discussed in the present paper is different. The objective of this
work is to rigorously prove that, in our setting, the patterns are
created by the combined action of the slow diffusion and the
switching process. Furthermore, the proof allows us to determine
the time scale of pattern formation and some of the pattern
parameters.
Our results show that there is no single winner in
the evolutionary game of the two competing phenotypes, but rather
both phenotypes are present in the limit distribution
of bacteria over the interval of available threshold values. This
supports the idea that the bet-hedging strategy can be optimal
for increasing fitness. It is shown that the complexity
of the pattern in the limit distribution increases with
decreasing diffusion rate in agreement with the numerical results
presented in \cite{1, 2}.
Further, we show that a
steady state solution achieved in the long time limit depends on
the initial state. This effect can be interpreted as a manifestation
of the long term memory on the level of populations and agrees
with earlier findings suggesting that memory in the switching
strategy of individual bacteria can grant fitness in varying environments \cite{FD}.
Our results suggest that the attractor of the system is a connected continual set of stationary patterns.

The paper is organized as follows. In
Section~\ref{secModelDescription}, we introduce the model: a
system of reaction-diffusion equations with spatially distributed
relays, and recall the well-posedness property from~\cite{2}. In
Section~\ref{secMainResults}, we formulate our main results:
pattern formation on two different time scales, depending on the
initial data. In Section~\ref{secFundamSol}, we establish some
properties of solutions of the heat equation with initial data
close to the delta function. These properties play an essential
role in Sections~\ref{proof} and~\ref{proof2}, where we prove the
two main results about pattern formation.

\section{Model description}\label{secModelDescription}

\subsection{Model assumptions}
In this paper, we consider a class of models, which attempt to
account for a number of phenomena listed above, namely (a)
switching of bacteria between two phenotypes (states) in response
to variations of environmental conditions; (b) hysteretic
switching strategy (switching rules) associated with bi-stability
of phenotype states; (c) heterogeneity of the population in the
form of a distribution of switching thresholds; (d) bet-hedging in
the form of diffusion between subpopulations characterized by
different bi-stability ranges; and, (e) competition for nutrients.
The resulting model is a reaction-diffusion system including, as
reaction terms, discontinuous hysteresis relay operators and the
integral of those. This integral can be interpreted as the
Preisach operator \cite{SH} with a time dependent density (the
density is a component of the solution describing the varying
distribution of bacteria). The main objective of this paper is to
prove that fitness, competition and diffusion can act together to
select a nontrivial distribution of phenotypes (states) over the
population of thresholds.

We assume that each of the two phenotypes, denoted by $1$ and
$-1$, consumes a different type of nutrient (for example, one
consumes lactose and the other glucose). The amount of nutrient
available for phenotype $i$ at the moment $t$ is denoted by
$f_i(t)$ where $i=\pm 1$. The model is based on the following
assumptions (see~\cite{2} for further discussion).

\begin{enumerate}
\item Each bacterium changes phenotype in response to the
variations of the variable $w=f_1/(f_1+f_{-1})-1/2$, which
measures the deviation of the relative concentration of the first
nutrient from the value 1/2 in the mixture of the two nutrients.

\item The input $w=w(t)$ is mapped to the (binary)   phenotype
(state) of a bacterium $r(t)={\mathcal R^x}(w)(t)$, where
${\mathcal R^x}$ is the non-ideal relay operator   with symmetric
switching thresholds $x$, $y=-x$ with $x>0$; see Fig.~\ref{hyst}
and the rigorous definition~\eqref{nonideal} in
Section~\ref{secSetting}.

\item The population includes bacteria with different bi-stability
ranges $(-x,x)$, where the threshold value $x$ varies over an
interval $[\underline{x}, \overline{x}]\subset (0,1/2)$. We will
denote by $u(x,t)$ the density of the biomass of bacteria with
given switching thresholds $\pm x$ at a moment $t$.

\item There is a diffusion process acting on the density $u$.

\item At any particular time moment $t$, for any given $x$, all
the bacteria with the switching threshold values $\pm x$ are in
the same state (phenotype).
   That is, $u(x,t)$ is the total density of bacteria with the threshold $x$ at the moment $t$ and they are all in the same state.
This means that when a bacterium with a threshold $x'$
sporadically changes its threshold to a different value $x$, it
simultaneously copies the state from other bacteria which have the
threshold~$x$. In particular, this may require a bacterium to
change the state when its threshold changes.
\end{enumerate}

With these assumptions, we obtain the following model of the
evolution of bacteria and nutrients,
\begin{equation}\label{00}
\left\{
\begin{aligned}
& {u}_t=  D u_{xx} + \frac12  (1+{\mathcal R}^x(w))\, u
f_1+\frac12(1-{\mathcal R}^x(w))\, u f_{-1},\\
& \dot{f}_{1} =-\frac12
f_1\int_{\underline{x}}^{\overline{x}}(1+{\mathcal R}^x(w))\,
u\,dx,\\
& \dot{f}_{-1} =-\frac12
f_{-1}\int_{\underline{x}}^{\overline{x}}(1-{\mathcal R}^x(w))\,
u\,dx,
\end{aligned}
\right.
\end{equation}
where $u_t$ and $u_{xx}$ are the derivatives of the population
density $u$;\ $D>0$ is the diffusion coefficient; dot denotes the
derivative with respect to time; and all the non-ideal relays
${\mathcal R^x}$, $x\in[\underline{x},\overline{x}]$, have the
same input $w=f_1/(f_1+f_{-1})-1/2$. We additionally assume the
growth rate $\frac12 (1+ i {\mathcal R}^x(w))\, u f_i $ based on
the mass action law for bacteria in the phenotype $i=\pm 1$.  The
rate of the consumption of nutrient in the equation for
$f_i=f_i(t)$ is proportional to the total biomass of bacteria in
the phenotype $i$, hence the integral; $\underline{x}$ and
$\overline{x}$ are the lower and upper bounds on available
threshold values, respectively.

We assume that a certain amount of nutrients is available at the
initial moment; the nutrients are not supplied after that moment.
 We assume the Neumann boundary conditions for $u$, that is no flux
of the population density $u$ through the lower and upper bounds
of available threshold values.

\subsection{Rigorous model setting}\label{secSetting}
Throughout the paper,  we assume  that
$x\in[\ux,\ox]\subset(0,1/2)$.

We begin with a rigorous definition of the hysteresis operator
$\cR^x$ ({\em non-ideal relay}) with fixed thresholds $\pm x$.
This  operator takes a continuous function $w=w(t)$ defined on an
interval $t\in[0,T)$ to the binary function $r=\cR^x(w)$ of time
defined on the same interval, which is given by
\begin{equation}\label{nonideal}
 {\mathcal R}^x (w)(t) = \left\{
\begin{array}{rll} -1 & {\rm if} & w(\tau)\le -x \ {\rm for\ some} \ \tau\in [0,t] \\
 && {\rm and} \ w(s)<x \ {\rm for\ all} \ s\in [\tau,t],
 \\1& {\rm if} &  w(\tau)\ge x \ {\rm for\ some} \ \tau\in [0,t] \\
 && {\rm and} \ w(s)>-x \ {\rm for\ all} \ s\in [\tau,t],\\
 r_0 & {\rm if} & -x < w(\tau) < x \ {\rm for\ all} \ \tau\in
 [0,t],
 \end{array}
\right.
\end{equation}
where $r_0$ is either $1$ or $-1$ (initial state of the non-ideal
relay $\cR^x$). Since $r_0$ may take different values for
different $x$, we write $r_0=r_0(x)$. The function $r_0=r_0(x)$ of
$x\in[\underline{x},\overline{x}]$ taking values $\pm 1$ is called
the {\em initial configuration} of the non-ideal relays. In what
follows, we  do not explicitly indicate the dependence of the
operator ${\mathcal R}^x$ on $r_0(x)$.

In this paper, we assume that $r_0(x)$ is {\it simple}, which
means the following. There is a partition $\underline{x}=\bar
x_0<\bar x_1<\cdots<\bar x_{N_0}=\overline{x}$ of the interval
$[\underline{x},\overline{x}]$ such that the function $r_0(x)$,
which satisfies $|r_0(x)|=1$ for all $\underline{x}\le x\le
\overline{x}$, is constant on each interval $(\bar x_{k-1},\bar
x_k]$ and has different signs on any two adjacent intervals:
\begin{equation}\label{simple}
\begin{array}{l}
r_0(x)=r_0(\bar x_k),\qquad x\in(\bar x_{k-1},\bar x_k],\ k=1,\ldots,N_0,\\
r_0(\bar x_{k-1})r_0(\bar x_k)=-1,\ k=2,\ldots,N_0,
\end{array}
\end{equation}
where the second relation holds if $N_0\ge2$.

We define the {\em distributed relay operator} $\cR(w)$ taking
functions $w=w(t)$ to functions $r=r(x,t)$ by
$$
r(x,t)=\cR(w)(x,t):={\mathcal R}^x(w)(t).
$$
The function $r(\cdot,t)$ will be referred to as the configuration
(state) of the distributed relay operator at the moment $t$.

We set
\begin{equation}\label{eqNP}
\uU(t)= \int_{\underline{ x}}^{\overline{ x}} u( x,t)\,d x,\quad
\cP(u,w)(t)=\int_{\underline{ x}}^{\overline{ x}} u( x,t)
\cR^x(w)(t)\,d x.
\end{equation}
Here the first integral is the total mass of bacteria and $\cP$ is
the so-called {\em Preisach operator} \cite{SH} with the time
dependent density function $u$; differential equations with the
Preisach operator have been studied, for example, in \cite{h1, h2,
h3, h4, h5, h6, h7, h8, h9, h10}. Further, we replace the unknown
functions $f_1$ and $f_{-1}$ in system \eqref{00} by
$v=f_1+f_{-1}$ (total mass of the two nutrients) and
$w=f_1/(f_1+f_{-1})-1/2$ (deviation of the relative concentration
of the first nutrient from the value~1/2). The resulting system
takes the form
\begin{equation}\label{eqBacteriaGeneral}\left\{
\begin{aligned}
& u_t = D u_{xx} +\left(\frac12 + w  \cR(w)\right)u v,\\
&\dot v=-\left(\frac{\underline{U}}{2}  +w
\cP(u,w)\right)v,\\
&\dot w=- \left(\dfrac{1}{2}+w\right)\left(\dfrac{1}{2}-w\right)
\cP(u,w),
\end{aligned}\right.
\end{equation}
where we assume the Neumann boundary conditions
\begin{equation}
\label{eqBC}  u_x|_{ x=\ux} = u_x|_{ x=\ox}=0
\end{equation}
and the initial conditions
\begin{equation}
\label{eqIC}
\begin{gathered}
u( x,0)=u_0( x),\quad v(0)=v_0,\quad w(0)=w_0,\quad
r(x,0)=r_0(x),\\
u_0( x)\ge0,\quad \int_\ux^\ox u_0(x)\,dx =1,\quad v_0\ge 0, \quad
|w_0|\le \ox,\quad r_0(x)\ \text{is simple}.
\end{gathered}
\end{equation}

\subsection{Well-posedness}\label{secWellPosedness}
The problem  \eqref{eqBacteriaGeneral}--\eqref{eqIC}, which
contains a discontinuous distributed relay operator, was shown  in
\cite{2} to be well posed. We briefly summarize this result before
proceeding with the analysis of long time behavior.

Set  $Q_T=(\ux,\ox)\times(0,T)$ for $T>0$. We will use the
standard Lebesgue spaces $L_2(Q_T)$ and $L_ 2=L_ 2(\ux,\ox)$; the
Sobolev spaces $W_ 2^k=W_ 2^k(\ux,\ox)$, $k\in\bbN$; the
anisotropic Sobolev space $W_2^{2,1}(Q_T)$  with the norm
$$
\|u\|_{W_2^{2,1}(Q_T)}=\left(\int\limits_0^T
\|u(\cdot,t)\|_{W_2^2}^2\, dt+ \int\limits_0^T
\|u_t(\cdot,t)\|_{L_2}^2\, dt\right)^{1/2};
$$
and, the space
$$
\cW(Q_T)=W_2^{2,1}(Q_T)\times C^1[0,T]\times C^1[0,T].
$$

\begin{definition}
Assume that $(u_0,v_0,w_0) \in W_2^1\times\bbR^2$. We say that
$(u,v,w)$ is a ({\em global})  {\em solution} of
problem~\eqref{eqBacteriaGeneral}--\eqref{eqIC}  if, for any
$T>0$,
 $(u,v,w)\in
\cW(Q_T)$, $r(\cdot,t)$ is a continuous $L_2$-valued function for
$t\ge0$, and relations~\eqref{eqBacteriaGeneral}--\eqref{eqIC}
hold in the corresponding function spaces.
\end{definition}


The following result was proved in \cite{2}.

\begin{theorem}\label{thWellPosed}
If  $(u_0,v_0,w_0) \in W_2^1\times\bbR^2$, then
\begin{enumerate}
\item problem~$\eqref{eqBacteriaGeneral}$--$\eqref{eqIC}$ has a
unique solution $(u,v,w)$$;$

\item the state $r(\cdot,t)=\cR(w)(\cdot,t)$ of the distributed
relay operator is {simple} for all $t\ge 0$$;$

\item we have
\begin{gather} \dot \uU(t)\ge0,\quad \uU(t)\to
1+v_0,\quad
 0\le u(\cdot,t)\to\dfrac{1+v_0}{\ox-\ux}\ \ \text{in}\ \
 C[\ux,\ox], \label{limit''}\\
\dot v(t)\le 0,\qquad v(t)\le v_0e^{-\mu t} \quad {\rm with} \quad \mu=1/2-\overline{x}>0, \label{limit'} \\
 |w(t)|\le\ox<1/2, \label{exist}
\end{gather}
where   the convergence takes place as $t\to\infty$.
\end{enumerate}
\end{theorem}





The   behavior given by~\eqref{limit''}--\eqref{exist} is to be
expected. Indeed, as we assume no supply of nutrients after the
initial moment, the total amount of nutrients $v(t)$ converges to
zero. When the density of nutrients vanishes as a result of
consumption by bacteria, the equation for the density $u$
approaches the homogeneous heat equation with zero flux boundary
conditions, which explains why the density of bacteria $u(x,t)$
converges to a uniform distribution over the interval
$[\underline{x},\overline{x}]$ as a result of the diffusion.

\section{Main results}\label{secMainResults}


The main observation we make in this paper is that for some
initial data \eqref{eqIC}, the limit distribution $r_*(x)$ has
many sign changes provided that the rate  $D$ of the diffusion
process is small enough. In particular, if  $r_0(x)\equiv 1$ or
$r_0(x)\equiv -1$, we interpret the transition of system
\eqref{eqBacteriaGeneral} from a uniform state $r_0(x)$ to the
sign changing state $r_*(x)$ as the formation of a pattern in the
distribution of phenotypes over the set of thresholds in bacteria.

Below, we consider the initial density satisfying the additional
relation
\begin{equation}\label{data2}
u_0(x)\le\varepsilon,\quad  x\in
(\underline{x},\overline{x}-\varepsilon)
\end{equation}
for some $\varepsilon>0$.

\begin{theorem}\label{t1} For any $N\in\bbN$ and
$v_0\ge 0$, there exist $D_*=D_*(N,v_0)>0$,
$\varepsilon_*=\varepsilon_*(N,v_0)>0$, and   sequences
$\{\tau_i\}_{i=1}^N$ and $\{x_i\}_{i=1}^{N}$ such that
$$
\begin{gathered}
0<\tau_N  <\tau_{N-1}<\dots<\tau_2< \tau_1,\\
\ux<x_1<x_2<\dots<x_N<\ox,
\end{gathered}
$$
and the following is true. If $D\in(0,D_*)$ and $u_0(x)$
satisfies~\eqref{data2} with $\varepsilon\in(0,\varepsilon_*)$,
then the state function $r(\cdot,t)$ of the distributed relay
operator has at least $N-i+1$ sign changes for $x\in [x_{i},\ox)$
and for all $t\ge\tau_{i}/D$.
\end{theorem}

Rephrasing the statement of this theorem, if the initial density
 $u_0(x)$ of bacteria approaches the delta function $\delta(x-\overline{x})$
and the diffusion $D$ goes to zero, then the number of sign
changes of the limit distribution $r_*(x)$ tends to infinity.


Theorem \ref{t1} establishes the pattern formation on the time
scale of order~$D^{-1}$. Under the additional assumption that the
amount of nutrients in the system is small, we show that the
pattern is formed on the time scale of order~$1$.

\begin{theorem}\label{tt1} There exist sequences $\{t_i\}_{i=1}^\infty$ and
$\{y_i\}_{i=1}^{\infty}$ satisfying
$$
\begin{gathered}
0< t_1< t_2<\dots,\qquad 0 <y_1<y_2<\dots
\end{gathered}
$$
such that, for any $N\in\bbN$, there are $D_*=D_*(N)>0$ and
$v_*=v_*(N)>0$ with the following property.  For any
$D\in(0,D_*)$, one can find $\varepsilon_*=\varepsilon_*(D)>0$
such that if $u_0(x)$ satisfies~\eqref{data2} with
$\varepsilon\in(0,\varepsilon_*)$ and $v_0\in[0, v_*)$, then the
state function $r(\cdot,t)$ of the distributed relay operator has
at least $i$ sign changes for $x\in [x_{i},\ox)$ and for all $t\ge
t_i$, where $x_i=\ox-D^{1/2}y_i$ and $i=1,\dots,N$.
\end{theorem}

\section{Fundamental solutions of the heat
equation}\label{secFundamSol}

\subsection{Setting of auxiliary problems for the heat equation}

In this section, we prove some auxiliary results that allow us to
compare solutions of the heat equation on finite intervals with
the fundamental solution on the infinite interval. We will use two
time scales given by the variable $t$ for systems with arbitrary
diffusion $D$ and by the variable $\tau$ for systems with
diffusion equal to $1$. Later, we will relate these time scales
via $\tau=Dt$.

{\bf Problem 1}:
\begin{equation}\label{eqpsi}
 \left\{
\begin{aligned}
&\psi_\tau=\psi_{xx},\quad x\in(-\infty,\overline{x}),\ \tau>0,\\
&\psi_x|_{x=\ox}=0,\quad \psi|_{\tau=0}=\delta(x-\overline{x}),
\end{aligned}\right.
\end{equation}
where $\delta(\cdot)$ is the Dirac delta function. Its ({\em
fundamental\/}) solution is well known to be given by
\begin{equation}\label{eqpsixtau}
\psi(x,\tau)=\frac{1}{\sqrt{\pi \tau}}\,
\exp\left(-\frac{(\overline{x}-x)^2}{4\tau}\right).
\end{equation}

{\bf Problem 2}:
\begin{equation}\label{eqphif}
\left\{
\begin{aligned}
&\varphi_\tau=\varphi_{xx}, \quad x\in(\ux,\overline{x}),\
\tau>0,\\
&\varphi_x|_{x=\ox}=\varphi_x|_{x=\ux}=0,\quad
\varphi|_{\tau=0}=\delta(x-\overline{x}).
\end{aligned}\right.
\end{equation}
Formally solving this problem by the Fourier method, we obtain the
solution
\begin{equation}\label{eqphi}
\varphi(x,\tau)=\frac{1}{\overline{x}-\underline{x}}+\frac{2}{\overline{x}-\underline{x}}\sum_{n=1}^\infty
\exp\left(-\left(\frac{n\pi}{\overline{x}-\underline{x}}\right)^2\tau\right)\cos
\left(\frac{n \pi  (\overline{x}-x)}{\overline{x}-\underline{x}}
\right),
\end{equation}
which we refer to as the {\it fundamental solution of
problem\/}~\eqref{eqphif}. Obviously, for any $\hat\tau>0$, the
series converges uniformly in the region
\begin{equation}\label{eqRegionuxoxtau0}
\{(x,\tau): x\in[\ux,\ox],\ \tau\ge\hat \tau\}
\end{equation}
and defines there a bounded infinitely differentiable function
$\varphi(x,\tau)$, which satisfies the heat equation and the
Neumann boundary conditions in the classical sense. Using
Lemma~\ref{lphi-psi} below, one can easily see that the initial
condition $\varphi|_{\tau=0}=\delta(x-\overline{x})$ is satisfied
in the following sense: for any $g\in C[\ux,\ox]$,
$$
\int_\ux^\ox \varphi(x,\tau)g(x)\,dx\to g(\ox)\quad\text{as }
\tau\to 0.
$$

{\bf Problem 3}:
\begin{equation}\label{equD}
\left\{
\begin{aligned}
&\omega_t=D\omega_{xx}, \quad x\in(\ux,\overline{x}),\
t>t_0,\\
&\omega_x|_{x=\ox}=\omega_x|_{x=\ux}=0,\quad
\omega|_{t=t_0}=\omega_0(x),
\end{aligned}\right.
\end{equation}
where $t_0\ge 0$ and $\omega_0\in W_2^1(\ux,\ox)$.  Set
$$
\Omega_0=\int_{\ux}^{\ox}\omega_0(x)\,dx.
$$
We assume that
\begin{equation}\label{data2uu}
0\le\dfrac{\omega_0(x)}{\Omega_0}\le \varepsilon_0,\quad x\in
[\underline{x},\overline{x}-\delta_0],
\end{equation}
for some $\varepsilon_0,\delta_0>0$. Using the Fourier method, we
represent the solution of problem~\eqref{equD} as the series
\begin{equation}\label{eqomega}
\omega(x,t;\varepsilon_0,\delta_0,D,t_0)=\frac{\Omega_0}{\overline{x}-\underline{x}}+\frac{2}{\overline{x}-\underline{x}}\sum_{n=1}^\infty
  A_n
\exp\left(-\left(\frac{n\pi}{\overline{x}-\underline{x}}\right)^2
D(t-t_0)\right)\cos \left(\frac{n \pi
(\overline{x}-x)}{\overline{x}-\underline{x}} \right)
\end{equation}
with the coefficients
\begin{equation}\label{eqAn}
A_n=\int_{\underline{x}}^{\overline{x}} \omega_0(x) \cos
\left(\frac{n \pi (\overline{x}-x)}{\overline{x}-\underline{x}}
\right)\,dx.
\end{equation}
We explicitly indicate the dependence of $\omega$ on the
parameters $\varepsilon_0,\delta_0,D,t_0$. We will see later on,
in Lemma~\ref{lphi-omega1}, that these parameters (rather than the
explicit form of the initial data $\omega_0(x)$) determine the
closeness between $\omega$ in Problem 3 and $\varphi$ in Problem
2.

Note that the series in~\eqref{eqomega} converges in
$W_2^{2,1}((\ux,\ox)\times(t_0,T))$ and hence in
$C([\ux,\ox]\times[t_0,T])$ for any $T>t_0$. It also
  converges uniformly in the region~\eqref{eqRegionuxoxtau0} and defines there a bounded
 function $\omega(x,t;D,t_0)$ that is infinitely differentiable with respect to $(x,t)$.

\subsection{Comparing Problems 1, 2, and 3}

First we compare Problems 1 and 2.

\begin{lemma}\label{lphi-psi} For any $\theta>0$,
\begin{equation}\label{lem}
\|\varphi(\cdot,\tau)-\psi(\cdot,\tau)\|_{C[\underline{x},\overline{x}]}\le
c(\theta) \exp\left(-\frac{(\ox-\ux)^2}{{5\tau}}\right),\quad
\tau\in(0,\theta],
\end{equation}
 where $c(\theta) >0$   does not depend on $\tau\in(0,\theta]$.
\end{lemma}
\proof Let us expand $\psi(x,\tau)$ given by~\eqref{eqpsixtau}
into the Fourier series
\begin{equation}\label{eqPsiFourier}
\psi(x,\tau)=\frac{1}{\overline{x}-\underline{x}}\sum_{n=0}^\infty
\psi_n(\tau)\cos\left(\frac{n \pi
(\overline{x}-x)}{\overline{x}-\underline{x}}\right)
\end{equation}
on the interval $\underline{x}\le x\le \overline{x}$. Here
\begin{equation}\label{eqPsi0}
\psi_0(\tau)=\dfrac{1}{\sqrt{\pi\tau}}\int\limits_0^{\overline{x}-\underline{x}}
\exp\left(-\frac{x^2}{4\tau}\right)\, dx
=1-\dfrac{2}{\sqrt{\pi}}\int\limits_{\frac{\overline{x}-\underline{x}}{2\sqrt{\tau}}}^\infty
\exp\left(-{x^2}\right)\,dx
=1-o\left(\exp\left(-\frac{(\overline{x}-\underline{x})^2}{{4\tau}}\right)\right),
\end{equation}
where $o(\cdot)$ is taken as $\tau\to 0$, and, for $n=1,2\dots$,
\begin{equation}\label{eqPsin1}
\psi_n(\tau)=\frac{4}{\sqrt{\pi}}
\int\limits_0^{\frac{\overline{x}-\underline{x}}{2\sqrt{\tau}} }
\exp\left(-x^2\right)\cos \left( \frac{2n \pi x\sqrt{\tau}}
{\overline{x}-\underline{x}} \right)\, dx.
\end{equation}
It follows from the identity
$$
\dfrac{1}{\sqrt{2\pi}}\int\limits_{-\infty}^{\infty}\exp\left({-\frac{x^2}{2}}\right)
\exp\left({i\xi x}\right)\,dx=\exp\left({-\frac{\xi^2}{2}}\right)
$$
  that
\begin{equation}\label{eqPsin2}
\frac{4}{\sqrt{\pi}} \int\limits_0^\infty
\exp\left({-x^2}\right)\cos \left( \frac{2n \pi x\sqrt{\tau}}
{\overline{x}-\underline{x}} \right)\, dx= 2
\exp\left({-\left(\frac{n\pi}{\overline{x}-\underline{x}}\right)^2\tau}\right).
\end{equation}
Further, using integration by parts and the estimate
$\int\limits_y^\infty
(1-2x^2)\exp\left({-x^2}\right)\,dx=O\left(y\,
\exp\left({-y^2}\right)\right)$ as $y\to\infty$, we obtain
\begin{equation}\label{eqPsin3}
\begin{aligned}
& \frac{4}{\sqrt{\pi}}
\int\limits_{\frac{\overline{x}-\underline{x}}{2\sqrt{\tau}}
}^\infty \exp\left({-x^2}\right)\cos \left( \frac{2n \pi
x\sqrt{\tau}} {\overline{x}-\underline{x}} \right)\, dx=
\dfrac{4(\ox-\ux)}{\pi^{3/2}n\sqrt{\tau}}\int\limits_{\frac{\overline{x}-\underline{x}}{2\sqrt{\tau}}
}^\infty x \exp\left({-x^2}\right)\sin \left( \frac{2n \pi
x\sqrt{\tau}}
{\overline{x}-\underline{x}} \right)\, dx\\
&\quad =\dfrac{2(\ox-\ux)^2}{\pi^{5/2}n^2
\tau}\left(\frac{\ox-\ux}{2\sqrt{\tau}}
\exp\left({-\frac{(\ox-\ux)^2}{{4\tau}}}\right)\cos(n\pi)+\int\limits_{\frac{\overline{x}-\underline{x}}{2\sqrt{\tau}}
}^\infty (1-2x^2)\exp\left({-x^2}\right)\cos \left( \frac{2n \pi
x\sqrt{\tau}}
{\overline{x}-\underline{x}} \right)\, dx\right)\\
&\quad
=\dfrac{1}{n^2\tau^{3/2}}O\left(\exp\left({-\frac{(\ox-\ux)^2}{{4\tau}}}\right)\right),
\end{aligned}
\end{equation}
where $O(\cdot)$ is taken as $\tau\to 0$ and does not depend on
$n$. Relations \eqref{eqPsin1}--\eqref{eqPsin3} imply
$$
\psi_n(\tau)=2
\exp\left({-\left(\frac{n\pi}{\overline{x}-\underline{x}}\right)^2\tau}\right)-\dfrac{1}{n^2}\,
O\left(\exp\left({-\frac{(\ox-\ux)^2}{{5\tau}}}\right)\right).
$$
Combining this formula with \eqref{eqphi}, \eqref{eqPsiFourier},
and \eqref{eqPsi0}, we conclude that
$$
|\varphi(x,\tau)-\psi(x,\tau)|=o\left(\exp\left({-\frac{(\ox-\ux)^2}{{4\tau}}}\right)\right)
+O\left(\exp\left({-\frac{(\ox-\ux)^2}{{5\tau}}}\right)\right)
\sum_{n=1}^\infty \dfrac{1}{n^2}.
$$
\endproof

Now we compare Problems 2 and 3.

\begin{lemma}\label{lphi-omega1}
Assume that \eqref{data2uu} holds with $\varepsilon_0\le 1$. Then,
for all $t>t_0$, we have
\begin{equation}\label{eqomega-phi1}
\begin{aligned}
&
\left\|\dfrac{\omega(\cdot,t;\varepsilon_0,\delta_0,D,t_0)}{\Omega_0}
-
\varphi(\cdot,D(t-t_0))\right\|_{C[\underline{x},\overline{x}]}\\
&\qquad \le C\left(\delta_0^2 (D(t-t_0))^{-3/2}+\varepsilon_0
(D(t-t_0))^{-1/2}\right),
\end{aligned}
\end{equation}
where $C>0$ does not depend on $t,\varepsilon_0,\delta_0,D,t_0$.
\end{lemma}
\proof Using~\eqref{data2uu}, one can easily see  that the
coefficients $A_n$ given by~\eqref{eqAn} satisfy
$$
\left|\dfrac{A_n}{\Omega_0}-1\right|\le
k_1(n^2\delta_0^2+\varepsilon_0),
$$
where $k_1>0$ does not depend on $n,\varepsilon_0,\delta_0$.
Therefore,\footnote{To prove the last inequality
in~\eqref{eqomega-phi1-proof}, we use the following. Since, for
any $\sigma>0$, the function $x\,\exp(-\sigma x)$ achieves its
maximum
 $\sigma^{-1}e^{-1}$ at $x=\sigma^{-1}$, we have
$$
\begin{aligned}
\sum\limits_{n=1}^\infty n^2 \exp(-n^2\sigma)& =
\sum\limits_{n=1}^{[\sigma^{-1/2}]} n^2
\exp(-n^2\sigma)+\sum\limits_{n=[\sigma^{-1/2}]+1}^\infty n^2
\exp(-n^2\sigma)\\
&\le e^{-1}\sigma^{-3/2}+ \int_{\sigma^{-1/2}}^\infty x^2
\exp(-x^2\sigma)\,dx
 =\left(e^{-1}+
\int_{1}^\infty y^2 \exp(-y^2)\,dy\right)\sigma^{-3/2}.
\end{aligned}
$$
This yields the desired estimate for the first sum
in~\eqref{eqomega-phi1-proof}. The estimate for the second sum
in~\eqref{eqomega-phi1-proof} follows from the inequality
$$
\sum\limits_{n=1}^\infty \exp(-n^2\sigma)\le\int_0^\infty
\exp(-x^2\sigma)\,dx=\dfrac{\sqrt{\pi}}{2}\sigma^{-1/2}.
$$
}
\begin{equation}\label{eqomega-phi1-proof}
\begin{aligned}
&\left\|\dfrac{\omega(\cdot,t;\varepsilon_0,\delta_0,D,t_0)}{\Omega_0}
- \varphi(\cdot,D(t-t_0))\right\|_{C[\underline{x},\overline{x}]}
\\
&\qquad \le k_1\left(\delta_0^2 \sum_{n=1}^\infty n^2
\exp\left({-\left(\frac{n\pi}{\overline{x}-\underline{x}}\right)^2
D(t-t_0)}\right)+\varepsilon_0\sum_{n=1}^\infty
\exp\left({-\left(\frac{n\pi}{\overline{x}-\underline{x}}\right)^2
D(t-t_0)}\right)\right)\\
&\qquad \le k_2\left(\delta_0^2 (D(t-t_0))^{-3/2}+\varepsilon_0
(D(t-t_0))^{-1/2}\right),
\end{aligned}
\end{equation}
where $k_2>0$ does not depend on $t>t_0$ and
$\varepsilon_0,\delta_0,D,t_0$.
\endproof

 \section{Proof of Theorem \ref{t1}}\label{proof}

\subsection{Observation moments}

We shall construct a sequence of time moments at each of which one
can observe a certain number of sign changes for the state
$r(\cdot,t)$ of the distributed relay operator $\cR(w)$. We call
these time moments {\it observation moments}. In this subsection,
we shall obtain estimates for the $u$-component of  the solution
$(u,v,w)$ of problem \eqref{eqBacteriaGeneral}--\eqref{eqIC}. They
provide information about the amount of bacteria on certain parts
of the interval $[\ux,\ox]$.  In the next subsection, we will use
this information in order to estimate from below the number of
sign changes for the state $r(\cdot,t):
[\underline{x},\overline{x}]\to \{-1,1\}$.

We set
\begin{equation}\label{eqUintu}
U(x,t)=\int_x^{\overline x} u(y,t)\,dy.
\end{equation}
 Then
$$
\uU(t)=U(\ux,t)=\int_\ux^{\overline x} u(y,t)\,dy;
$$
cf.~\eqref{eqNP}. The goal of this subsection is to prove the
following result.

\begin{lemma}\label{lUinequalities}
For any $N\in\bbN$, there exist $D_*=D_*(N,v_0)>0$,
$\varepsilon_*=\varepsilon_*(N,v_0)>0$, and positive sequences
$\{\tau_i\}_{i=1}^N$, $\{\chi_i\}_{i=1}^N$, and
$\{x_i\}_{i=1}^{N+1}$ such that
$$
\begin{gathered}
0<\tau_N-\chi_N<\tau_N<\tau_{N-1}-\chi_{N-1}<\tau_{N-1}<\dots<\tau_2<\tau_1-\chi_1<\tau_1,\\
\ux<x_1<x_2<\dots<x_{N+1}<\ox,
\end{gathered}
$$
and
\begin{equation}\label{2*}
\frac{U(x_{i+1},\tau/D)}{\uU(\tau/D)}<\frac12, \quad   \tau\ge
\tau_i-\chi_i.
\end{equation}
\begin{equation}\label{1*}
2\frac{U(x_i,\tau/D)}{\uU(\tau/D)}-2\frac{U(x_{i+1},\tau/D)}{\uU(\tau/D)}-1>\chi_i,
 \quad \tau\in[\tau_i-\chi_i,\tau_i],
\end{equation}
provided that $D\in(0,D_*)$ and $u_0(x)$ satisfies~\eqref{data2}
with $\varepsilon\in(0,\varepsilon_*)$.
\end{lemma}

We begin with estimates similar to~\eqref{2*} and \eqref{1*} for
the function
\begin{equation}\label{eqPhi}
\Phi(x,\tau)=\int_{x}^{\overline{x}} \varphi(y,\tau)\,dy.
\end{equation}

\begin{lemma}\label{lPhiinequalities}
There exist positive sequences $\{\tau_i\}_{i=1}^\infty$,
$\{\chi_i\}_{i=1}^\infty$, and $\{x_i\}_{i=1}^{\infty}$ such that
$\tau_i,\chi_i\to0$, $x_i\to\ox$ as $i\to\infty$,
\begin{equation}\label{eqtaux}
\begin{gathered}
0<\dots<\tau_i-\chi_i<\tau_i<\tau_{i-1}-\chi_{i-1}<\tau_{i-1}<\dots<\tau_2<\tau_1-\chi_1<\tau_1,\\
\ux<x_1<x_2<\dots<x_i<\dots<\ox,
\end{gathered}
\end{equation}
 and
\begin{equation}\label{2'}
\Phi(x_{i+1},\tau)<1/2 -\chi_i \quad {\rm for} \quad \tau\ge
\tau_i-\chi_i,
\end{equation}
\begin{equation}\label{1'}
2\Phi(x_i,\tau)-2\Phi(x_{i+1},\tau)-1>2\chi_i \quad {\rm for}
\quad \tau\in[\tau_i-\chi_i,\tau_i].
\end{equation}
\end{lemma}
\proof Note that
\begin{equation}\label{p1}
\Phi(x,\tau) \to 1 \quad {\rm as} \quad \tau\to 0  \quad \forall
x\in [\underline x,\overline{x});
\end{equation}
\begin{equation}\label{p2}
\sup _{\tau\ge \hat\tau} \Phi(x,\tau) \to 0 \quad {\rm as} \quad
x\to \overline{x}  \quad \forall\hat\tau>0.
\end{equation}
Indeed, set
$$
\Psi(x,\tau)=\int_{x}^{\overline{x}} \psi(y,\tau)\,dy.
$$
Then it is easy to see that $\Psi(x,\tau) \to 1$  as $\tau\to 0$
for any $x\in [\underline x,\overline{x})$. Combining this
with~\eqref{eqPhi} and Lemma~\ref{lphi-psi}, we obtain~\eqref{p1}.
Relation~\eqref{p2} follows from~\eqref{eqphi} and~\eqref{eqPhi}.

Now we construct sequences $\{\tau_i\}$ and $\{x_i\}$
satisfying~\eqref{eqtaux} and such that
\begin{equation}\label{2}
\Phi(x_{i+1},\tau)<1/2 \quad \forall \tau\ge \tau_i,
\end{equation}
\begin{equation}\label{1}
2\Phi(x_i,\tau_i)-2\Phi(x_{i+1},\tau_i)-1>0.
\end{equation}
Take an arbitrary $x_1\in(\underline{x},\overline{x})$ and, using
\eqref{p1}, choose   $\tau_1>0$ such that
$$
\Phi(x_1,\tau_1)>1/2.
$$
This estimate and relation \eqref{p2} ensure that inequality
\eqref{1} with $i=1$ holds for any $x_2$ which is sufficiently
close to $\overline{x}$. Simultaneously, relation~\eqref{p2}
ensures \eqref{2} for such $x_2$. Hence, we can find $x_2\in
(x_1,\overline{x})$ such that both \eqref{2} and \eqref{1} are
satisfied for $i=1$. Now, using \eqref{p1} again, we can choose a
sufficiently small $\tau_2<\tau_1$ such that
$$
\Phi(x_2,\tau_2)>1/2
$$
and continue constructing the sequences $x_i\to \overline{x}-,
\tau_i\to 0+$ by induction.

After the sequences $\{x_i\}$ and $\{\tau_i\}$ have been
constructed, we choose a  positive sequence $\chi_i$ such
that~\eqref{2'} and~\eqref{1'} hold and $
\tau_{i}<\tau_{i-1}-\chi_{i-1}. $
\endproof

In the next two lemmas, we show that
$\dfrac{U(x,\tau/D)}{\uU(\tau)}$ is close to $\Phi(x,\tau/D)$.
Together with Lemma~\ref{lPhiinequalities}, this will complete the
proof of Lemma~\ref{lUinequalities}.

We proceed in two steps. In  the first step
(Lemma~\ref{lsmallutau0D}), we show that one can choose $\tau_0>0$
not depending on $D$ such that the density $u$ at the time moment
$\tau_0/D$ remains close to the delta function $\delta(x-\ox)$. In
the second step (Lemma~\ref{lUclosePhi}), using the exponential
decay of the amount of nutrient~$v(t)$, we choose $D$ sufficiently
small such that $v(t)$ at the moment $\tau_0/D$ is close to zero
and show that, after this time moment, the density $u(x,t)$ can be
approximated by the fundamental solution $\varphi(x,t)$.

\begin{lemma}\label{lsmallutau0D}
For any   $\varepsilon_1>0$, there is a sufficiently small
$\theta_0=\theta_0(\varepsilon_1,v_0)>0$ with the following
property. Given $\tau_0\le\theta_0$, there exists
$\varepsilon_*=\varepsilon_*(\tau_0,v_0)>0$ such that
\begin{equation}\label{u2}
u(x,\tau_0/D)\le\varepsilon_1 \quad \text{for all}\
x\in[\underline{x},\overline{x}-\varepsilon_1],\ D>0,
\end{equation}
provided that the initial data $u_0(x)$ satisfies~\eqref{eqIC}
and~\eqref{data2} with $\varepsilon\in(0,\varepsilon_*)$.
\end{lemma}
\proof Set
$$
a(x,t)=\left(\frac12+w(t){\mathcal R}^x (w)(t)\right)v(t).
$$
Then the component $u$ of the solution $(u,v,w)$ solves the linear
problem
\begin{equation}\label{cv}
u_t=D u_{xx} + a(x,t)u,\quad u_x|_{x=\ux}=u_x|_{x=\ox}=0,\quad
u|_{t=0}=u_0(x).
\end{equation}
Combining the equality $|\cR^x(x)(t)|=1$ with
relations~\eqref{limit'} and~\eqref{exist}, we see that $0\le
a(x,t)\le v_0 e^{-\mu t}$. It follows from the comparison theorems
that
\begin{equation}\label{equlowerupper}
\uu(x,t)\le u(x,t)\le \ou(x,t),\quad
x\in[\underline{x},\overline{x}],\ t\ge0,
\end{equation}
where the lower and upper solutions $\uu$ and $\ou$ satisfy the
equations
\begin{equation}\label{probl}
\underline{u}_t=D \underline{u}_{xx},\qquad \overline{u}_t=D
\overline{u}_{xx} +v_0e^{-\mu t}\overline{u},
\end{equation}
respectively, with the same boundary and initial conditions as in
\eqref{cv}.

By Lemma~\ref{lphi-omega1}, for any $\tau_0>0$,
\begin{equation}\label{limi}
\|\underline {u}(\cdot,\tau_0/D) - \varphi
(\cdot,\tau_0)\|_{C[\underline{x},\overline{x}]}
  \le \varepsilon_0,
\end{equation}
where
\begin{equation}\label{eqeps0ofeps}
\varepsilon_0=C\left(\varepsilon^2 \tau_0^{-3/2}+\varepsilon
 \tau_0^{-1/2}\right)
\end{equation}
and $C>0$ does not depend on $\varepsilon,\tau_0,D$.
Using~\eqref{equlowerupper}, \eqref{limi}, and the equality
\begin{equation}\label{equuoue}
\overline{u}(x,t)=\underline{u}(x,t)
\exp\left({\frac{v_0(1-e^{-\mu t})}{\mu}}\right),
\end{equation}
we obtain
\begin{equation}\label{bound}
 u(x,\tau_0/D)\le
\exp\left({\frac{v_0}{\mu}}\right) (\varphi
(x,\tau_0)+\varepsilon_0),\quad x\in[\underline{x},\overline{x}].
\end{equation}

Using Lemma~\ref{lphi-psi}, it is easy to see that
$\sup\limits_{x\in[\ux,\hat x]}\varphi(x,\tau)\to 0$ as $\tau\to
0$ for any $\hat x\in[\ux,\ox)$. Therefore, for any
$\varepsilon_1>0$, there is a sufficiently small
$\theta_0=\theta_0(\varepsilon_1,v_0)>0$ such that
\begin{equation}\label{equ2_1}
\exp\left({\frac{v_0}{\mu}}\right)\varphi(x,\tau_0)\le\dfrac{\varepsilon_1}{2}\quad\text{for}\
x\in[\ux,\ox-\varepsilon_1],\ \tau_0\le\theta_0.
\end{equation}
Further, using~\eqref{eqeps0ofeps}, one can choose a sufficiently
small $\varepsilon_*=\varepsilon_*(\tau_0,v_0)$ such that
\begin{equation}\label{equ2_2}
\exp\left({
\frac{v_0}{\mu}}\right)\varepsilon_0\le\dfrac{\varepsilon_1}{2},
\end{equation}
provided that $\varepsilon\le\varepsilon_*$.
Combining~\eqref{bound}--\eqref{equ2_2} yields~\eqref{u2}.
\endproof

\begin{lemma}\label{lUclosePhi}
For any $N>0$ and $\varepsilon_2>0$, there are
$\varepsilon_*=\varepsilon_*(N,v_0,\varepsilon_2)>0$ and
$D_*=D_*(N,v_0,\varepsilon_2)>0$ such that
$$
\left\|\dfrac{U(\cdot,\tau/D)}{\uU(\tau/D)} -
\Phi(\cdot,\tau)\right\|_{C[\underline{x},\overline{x}]}\le
 \varepsilon_2\quad \forall \tau\ge\tau_N-\chi_N,
$$
where $\tau_N$ and $\chi_N$ are defined in
Lemma~\ref{lPhiinequalities}, provided that $D\in(0,D_*)$ and the
initial data $u_0(x)$ satisfies~\eqref{eqIC} and~\eqref{data2}
with $\varepsilon\in(0,\varepsilon_*)$.
\end{lemma}
\proof  Take $\tau_0$ in Lemma~\ref{lsmallutau0D} such that
$\tau_0\le(\tau_N-\chi_N)/2$ and repeat the same line of argument
as in the proof of Lemma~\ref{lsmallutau0D}, but now for the time
interval $[\tau_0/D,\infty)$ to obtain an analogue of estimates
\eqref{bound} for $t\ge (\tau_N-\chi_N)/D$. The lower and upper
solutions $\underline{u}_*$ and $\overline{u}_*$ in the
inequalities
\begin{equation}\label{equ*u}
\underline{u}_*(x,t)\le u(x,t)\le \overline{u}_*(x,t),\qquad
x\in[\underline{x},\overline{x}],\ t\ge\tau_0/D,
\end{equation}
are defined by the same problems \eqref{probl} (with the Neumann
boundary conditions), but with the initial conditions
$$
\underline{u}_*|_{t=\tau_0/D}=\overline{u}_*|_{t=\tau_0/D}={u}|_{t=\tau_0/D}.
$$
In particular, this implies that
\begin{equation}\label{estim}
\overline{u}_*(x,t)= \underline{u}_*(x,t) \exp\left({ v_0\frac{
e^{-\mu \tau_0/D}-e^{-\mu t}}{\mu} }\right)\le
\underline{u}_*(x,t) \exp\left({ \frac{v_0 e^{-\mu
\tau_0/D}}{\mu}}\right).
\end{equation}

Again applying Lemma~\ref{lphi-omega1} and using the inequalities
$\dot\uU(t)\ge 0$ and $\tau_0\le(\tau_N-\chi_N)/2$, we obtain
\begin{equation}\label{difer0}
\begin{aligned}
&\left\|\frac{\underline{u}_*(\cdot,\tau/D)}{\uU(\tau_0/D)}-\varphi
(\cdot,\tau-\tau_0)\right\|_{C[\underline{x},\overline{x}]}\\
&\qquad \le C_0\left(\varepsilon_1^2
(\tau_N-\chi_N)^{-3/2}+\varepsilon_1
 (\tau_N-\chi_N)^{-1/2}\right)\le C_1 \varepsilon_1
\end{aligned}
\end{equation}
for all $\tau\ge\tau_N-\chi_N$, where $C_0,C_1>0$ do not depend on
$\tau,\tau_0,D,\varepsilon,\varepsilon_1$ (but $C_1$ depends on
$N$).

Next, we take $\tau_0=\tau_0(\varepsilon_1)$  such that the
inequality
 \begin{equation}\label{difer}
\|\varphi (\cdot,\tau-\tau_0)-\varphi
(\cdot,\tau)\|_{C[\underline{x},\overline{x}]}\le
C_1\varepsilon_1\quad \forall \tau\ge\tau_N-\chi_N
\end{equation}
holds along with the inequality $\tau_0\le(\tau_N-\chi_N)/2$.
Then~\eqref{difer0} and~\eqref{difer} imply
$$
\left\|\frac{\underline{u}_*(\cdot,\tau/D)}{\uU(\tau_0/D)}-\varphi
(\cdot,\tau)\right\|_{C[\underline{x},\overline{x}]} \le 2C_1
\varepsilon_1\quad \forall \tau\ge\tau_N-\chi_N.
$$

Using the last inequality and estimates~\eqref{equ*u}
and~\eqref{estim}, we have
$$
\varphi (x,\tau)-2C_1\varepsilon_1 \le
\frac{u(x,\tau/D)}{\uU(\tau_0/D)} \le \exp\left({ \frac{
v_0e^{-\mu \tau_0/D} }{\mu}}\right) (\varphi
(x,\tau)+2C_1\varepsilon_1)
$$
for all $x\in [\underline{x},\overline{x}]$,
$\tau\ge\tau_N-\chi_N$. Integrating with respect to $x$ yields
\begin{equation}\label{comb2}
\Phi (x,\tau)-2C_1\varepsilon_1 \le
\frac{U(x,\tau/D)}{\uU(\tau_0/D)} \le \exp\left({ \frac{
v_0e^{-\mu \tau_0/D} }{\mu}}\right) (\Phi
(x,\tau)+2C_1\varepsilon_1)
\end{equation}
for all $x\in [\underline{x},\overline{x}]$,
$\tau\ge\tau_N-\chi_N$.

Since $e^{-\mu \tau_0/D}\to 0$, $\uU(\tau/D)\to 1+v_0$ as $D\to 0$
uniformly with respect to $\tau\ge\tau_0$, and $|\Phi(x,\tau)|\le
1$, it follow that there is $D_*=D_*(\tau_0,v_0,\varepsilon_1)$
such that
$$
\Phi (x,\tau)-3C_1\varepsilon_1 \le
\frac{U(x,\tau/D)}{\uU(\tau/D)} \le \Phi (x,\tau)+3
C_1\varepsilon_1,
$$
provided that $D\le D_*$. Setting
$\varepsilon_1=\varepsilon_2/(3C_1)$, we complete the proof.
\endproof

Now Lemma~\ref{lUinequalities} follows from
Lemmas~\ref{lPhiinequalities} and~\ref{lUclosePhi}.

\subsection{Proof of Theorem~\ref{t1}}\label{subsecProofTheoremt1}

\subsubsection{Terminology}

For the rest of the proof of the theorem we adopt the following
terminology. Points of the interval $\underline{x}< x\le
\overline{x}$ are colored in {\em white} (the relay is in
state~$-1$) and {\em black} (the relay is in state 1). Coloring
evolves in time. A point $\bar x_j=\bar x_j(t)$ which separates an
interval of white from an interval of black will be called a {\em
front}; cf.~\eqref{simple}. The total number of fronts can vary,
but, by Theorem~\ref{thWellPosed}, remains finite at all times
(equivalently, the state of the distributed relay operator remains
simple at all times). A front can either stay (a {\em steady}
front) or move right. That is, any $\bar x_j(t)$ is a
non-decreasing function on the time interval of its existence. A
front disappears if it hits another front or is hit by another
front. A front is called {\em immortal} if it cannot disappear.

The definition of relays implies the following dynamics of fronts
(equivalently, dynamics of the state $r=r(x,t)$ of the distributed
relay operator~${\mathcal R}$) in response to variations of the
input $w$. If at some moment $t$ the fronts create a partition
$\underline{x}<\bar x_1(t)<\cdots< \bar x_{\bar N(t)-1}<\bar
x_{\bar N(t)}=\overline{x}$ of the segment
$[\underline{x},\overline{x}]$ into alternating black and white
intervals, then, necessarily, $|w(t)|\le \bar x_1(t)$. The fronts
$\bar x_2,\ldots,\bar x_{\bar N(t)}$ cannot move, they are steady;
the most left front $\bar x_1(t)$ can move only if $|w(t)|= \bar
x_1(t)$, otherwise (i.e., if $|w(t)|<\bar x_1(t)$) it is also
steady.
When the front $\bar x_1(t)$ hits the front $\bar x_2(t)$, these
two fronts disappear. A new front is born when the input $w(t)$
leaves the interval $[-\underline{x},\underline{x}]$ for the first
time, and then each time the input leaves this interval after
having entered it from the other end.

\subsubsection{Collision moments}

An important moment for our consideration is when the front $\bar
x_1$ hits the front $\bar x_2$. At such a collision moment $t$ the
fronts $\bar x_1$ and $\bar x_2$ disappear and the front $\bar
x_3(t-)$ becomes the most left front, i.e., $ \bar x_1(t+)=\bar
x_3(t-)$.

\begin{lemma}\label{lfla}
Assume that the front $\bar x_1$ hits the front $\bar x_2$ at a
moment $t$. Then
\begin{equation}\label{fla}
\frac{U(\bar x_1(t+),t)}{\uU(t)}\ge\frac12.
\end{equation}
\end{lemma}
\proof Note that $w(t)$ and $\dot w(t)$ are continuous (see
Sec.~\ref{secWellPosedness}) and, at the collision moment $t$,
$w(t)\dot w(t)\ge 0$. Substituting the right-hand side of the last
equation of system \eqref{eqBacteriaGeneral} into this relation
and taking into account that $1/2\pm w(t)\ge\mu>0$ at all times,
we obtain
\begin{equation}\label{fla'}
0\ge w (t){\mathcal P}(u,w)(t)=w(t)\int_{\underline
x}^{\overline{x}} u(x,t){\mathcal R}^x (w)(t)\,dx.
\end{equation}

Further, for $s$ slightly less than the collision moment $t$, the
intervals $(\ux,\bar x_1(s))$ and $(\bar x_2(s), \bar x_3(s))$ had
the same color:
$$
{\mathcal R}^x (w)(s)={\rm sign}\, w(s), \quad x\in (\ux,\bar
x_1(s))\cup(\bar x_2(s), \bar x_3(s)),
$$
while the interval $(\bar x_1(s), \bar x_2(s))$ had the opposite
color:
$$
{\mathcal R}^x (w)(s)=-{\rm sign}\, w(s), \quad x\in  (\bar
x_1(s), \bar x_2(s)).
$$
At the collision moment $t$, the front $\bar x_3(t-)$ becomes the
most left front. Hence, we renumber the fronts. In particular, we
have $ \bar x_1(t+)=\bar x_3(t-)$, and now the whole interval
$(\ux,\bar x_1(t+))$ has the same color:
$$
{\mathcal R}^x (w)(t)={\rm sign}\, w(t), \quad x\in (\ux,  \bar
x_1(t+)).
$$
Combining this with \eqref{fla'}, we obtain
$$
0\ge \int_{\underline x}^{\bar x_1(t+)} u(x,t)\,dx+{\rm sign}\,
w(t)\int_{\bar x_1(t+)}^{\overline{x}}u(x,t){\mathcal R}^x
(w)(t)\,dx,
$$
where the first integral equals $\uU(t)-U(\bar x_1(t+),t)$ and the
absolute value of the second integral does not exceed $U(\bar
x_1(t+),t)$. That is, at the collision moment, inequality
\eqref{fla} must hold.
\endproof

\subsubsection{Completion of the proof of Theorem~\ref{t1}}  Fix an arbitrary
$N\in\bbN$. Then fix $D_*=D_*(N,v_0)>0$,
$\varepsilon_*=\varepsilon_*(N,v_0)>0$, and the sequences
$\{\tau_i\}_{i=1}^N$, $\{\chi_i\}_{i=1}^N$, and
$\{x_i\}_{i=1}^{N+1}$ as indicated in Lemma~\ref{lUinequalities}.
Then relations~\eqref{2*} and~\eqref{1*} hold for all
$D\in(0,D_*)$ and $u_0(x)$ satisfying~\eqref{data2} with
$\varepsilon\in(0,\varepsilon_*)$.

Now we will complete the proof by induction with respect to
$i=N+1,N,\dots,1$. Set $\tau_{N+1}=0$ and $t_{N+1}=0$.

{\it Basis.} The system has at least $0$ immortal fronts in the
interval $[x_{N+1}, \overline{x})$ at the moment $t_{N+1}=0$. This
statement is obviously true.

\begin{figure}[ht]
{\ \hfill\epsfxsize110mm\epsfbox{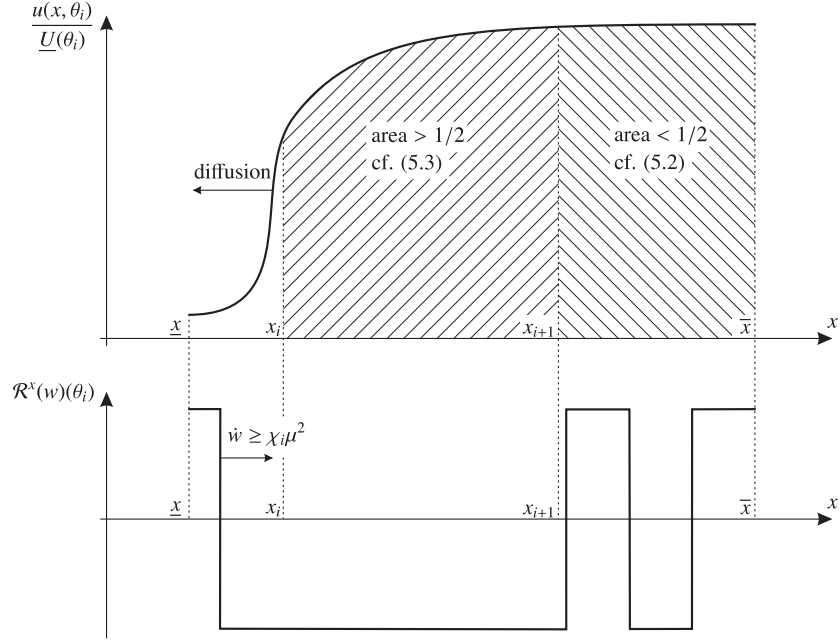}\hfill\ }
\caption{Spatial profiles of $ {u(x,t)}/{\uU(t)}$ and $\mathcal
R^x(w)(t)$ at the moment $t=\theta_i=(\tau_i-\chi_i)/D$, which
precedes $t_i=\tau_i/D$. The picture corresponds to the second
case in the inductive step: the most left front is not in the
interval $[x_i,x_{i+1})$, but, having a positive velocity,  will
enter $[x_i,x_{i+1})$ within the time interval
$t\in(\theta_i,t_i)$ and become immortal.} \label{FigSteps}
\end{figure}

{\it Inductive step: {\rm``}$N-i\ \Rightarrow\ N-i+1${\rm''} for
$1\le i\le N$.} We assume that the system has at least $N-i$
immortal fronts in the interval $[x_{i+1}, \overline{x})$ at the
moment $t_{i+1}=\tau_{i+1}/D$. Let us prove that it will have at
least $N-i+1$ immortal fronts in the interval $[x_i,\overline{x})$
after the moment $t_i=\tau_i/D$; see Fig.~\ref{FigSteps}.

First, we note that each front belonging to the interval
$[x_{i+1}, \overline{x})$ at the moment
$\theta_i=(\tau_i-\chi_i)/D$ is immortal after this moment.
Indeed, assuming the opposite, a collision of two fronts occurs at
some point of the interval $[x_{i+1},\overline{x}]$ at a moment
$t\ge\theta_i$. Therefore, by Lemma~\ref{lfla}, inequality
\eqref{fla} should be valid for a point $\bar x_1(t+)\in
[x_{i+1},\overline{x}]$  at the moment $t$. However, \eqref{fla}
contradicts \eqref{2*} as $U$ decreases in $x$.

Now two cases are possible.

In the first case, the interval $[x_i,x_{i+1})$ contains at least
one front at the moment $\theta_i=(\tau_i-\chi_i)/D$. Consider the
most right front~$F$ in this interval. We claim that the front $F$
is immortal. Indeed, $F$ cannot hit a front located to the right
of it because all the fronts to the right belong to
$[x_{i+1},\ox)$ and are  immortal by the above argument. Suppose,
at some moment $t\ge \theta_i$, the front $F$ is hit by another
front from the left. Then, at the collision moment, we have $\bar
x_1(t+)\ge x_{i+1}$  again. Hence, by Lemma~\ref{lfla},
inequality~\eqref{fla} should be valid. However, \eqref{fla}
contradicts \eqref{2*} because $U$ decreases in $x$.

In the second case, the interval $[x_i,x_{i+1})$  contains no
fronts at the moment $\theta_i=(\tau_i-\chi_i)/D$. Hence, all the
points of this interval are of the same color.
To be specific, assume without loss of generality that this
interval is white:
$$
{\mathcal R}^x(w)(\theta_i)=-1, \quad x\in [x_i,x_{i+1}).
$$
Hence, $w(\theta_i)<x_i$. Denote by $t_*$ the first moment when
$w(t)=x_i$ after the moment $\theta_i$, and set $t_*=\infty$ if
$w(t)$ never reaches $x_i$.

We claim that $t_*<t_i$. Indeed, during the time interval
$\theta_i\le t<t_*$, the interval $[x_i,x_{i+1})$ remains white.
Hence, using the last equation of
system~\eqref{eqBacteriaGeneral}, we obtain for  $\theta_i\le t <
\min\{t_*,t_i\}$
$$
\dot w =f(w)\left( \int_{x_i}^{x_{i+1}}
u(x,t)\,dx-\int_{\underline{x}}^{x_i} u(x,t){\mathcal
R}^x(w)\,dx-\int_{x_{i+1}}^{\overline{x}} u(x,t){\mathcal
R}^x(w)\,dx\right),
$$
where $f(w)=(1/2+w)(1/2-w)$. In this equation, the first integral
equals $U(x_i,t)-U(x_{i+1},t)$, while the sum of the absolute
values of the other two integrals does not exceed
$\uU(t)-U(x_i,t)+U(x_{i+1},t)$ (as $|{\mathcal R}^x(w)|=1$).
Therefore,
\begin{equation}\label{eqdotwineq1}
 \dot w(t)\ge
f(w)(2U(x_i,t)-2U(x_{i+1},t)-\uU(t)), \quad \theta_i\le t \le
\min\{t_*,t_i\}.
\end{equation}
Since estimate \eqref{1*} holds for all $t\ge \theta_i$, $f(w)\ge
\mu^2$, and $\uU(t)\ge 1$, it follows from~\eqref{eqdotwineq1}
that
\begin{equation}\label{eqdotwineq2}
\dot w(t)\ge \chi_i \mu^2, \quad \theta_i\le t \le
\min\{t_*,t_i\}.
\end{equation}
If $t_i\le t_*$, then inequality~\eqref{eqdotwineq2} implies
$$
1\ge w(t_i)-w(\theta_i)\ge \chi_i \mu^2
(t_i-\theta_i)=\chi_i^2\mu^2/D,
$$
where we have used the uniform estimate $|w|\le 1/2$. Hence,
decreasing $D_*$ if necessary, we obtain the contradiction.
Therefore, $\theta_i<t_*<t_i$. By definition of $t_*$, a front
enters the interval $[x_i, x_{i+1})$ from the lower end at the
moment~$t_*$. It is immortal by the same argument as in the first
case.

In both cases, we have proved that there will be at least $N-i+1$
immortal fronts in the interval $[x_i, \overline{x})$ after the
moment $t_*$ and hence after $t_i$. \qed

\section{Proof of Theorem \ref{tt1}}\label{proof2}
The proof is a modification of the above proof of Theorem
\ref{t1}.

\subsection{Observation moments}

The goal of this subsection is to prove an analogue of
Lemma~\ref{lUinequalities}, but keeping the observation moments of
order $1$ instead of $D^{-1}$ as $D\to 0$. Furthermore, the
distances between the immortal fronts will be of order $D^{1/2}$
instead of $1$.

 The main reason
of the time scale $O(D^{-1/2})$ in Lemma~\ref{lUinequalities} was
the necessity to wait long enough (namely, up to the time moment
$\tau_0/D$)  until the amount of nutrient $v(t)$ becomes small
enough. This is not necessary in the proof of Theorem \ref{tt1}
because we assume $v_0$ small from the very beginning.

As before, we set $\mu=1/2-\ox$.

\begin{lemma}\label{lUinequalitiesTheor2}
Fix $N\in\bbN$. Then there are $D_*=D_*(N)>0$, $v_*=v_*(N)>0$, and
sequences $\{t_i\}_{i=1}^N$, $\{\theta_i\}_{i=1}^N$, and
$\{y_i\}_{i=0}^{N}$ such that
\begin{equation}\label{eqty}
\begin{gathered}
0<\theta_1<t_1<\theta_2<t_2<\dots<\theta_N<t_N,\\
0=y_0<y_1<y_2<\dots<y_N,
\end{gathered}
\end{equation}
and the following is true. For any $D\in(0,D_*)$, one can find
$\varepsilon_*=\varepsilon_*(D)>0$ such that
\begin{equation}\label{5ol}
\frac{U(x_{i-1}, t)}{\uU(t)}< \frac{1}{2}, \quad t\ge \theta_i,
\end{equation}
\begin{equation}\label{4ol}
2\frac{U(x_i,  t)}{\uU(t)}-2\frac{U(x_{i-1},
t)}{\uU(t)}-1>\frac{1}{\mu^2(t_i-\theta_i)}, \quad t\in
[\theta_i,t_i],
\end{equation}
for $i=1,\dots,N$, where   $x_i=\ox-D^{1/2}y_i$, provided that
$u_0(x)$  satisfies~\eqref{data2} with
$\varepsilon\in(0,\varepsilon_*)$ and $v_0\in[0, v_*)$.
\end{lemma}

We begin with estimates similar to~\eqref{5ol} and \eqref{4ol} for
the function $\Phi$ given by~\eqref{eqPhi}.

\begin{lemma}\label{lPhiinequalitiesTheor2}
Fix $N\in\bbN$. Then there exist $D_*=D_*(N)>0$  and sequences
$\{t_i\}_{i=1}^N$, $\{\theta_i\}_{i=1}^N$, and $\{y_i\}_{i=0}^{N}$
satisfying~\eqref{eqty} and such that
\begin{equation}\label{3ol''}
\Phi(x_{i-1}, Dt)< \frac{1}{2}, \quad   t\ge \theta_i,
\end{equation}
\begin{equation}\label{3ol'}
2\Phi(x_i,D t)-2\Phi(x_{i-1},D t)-1>\frac{1}{\mu^2(t_i-\theta_i)},
\quad t\in [\theta_i,t_i],
\end{equation}
for $i=1,\dots,N$, where   $x_i=\ox-D^{1/2}y_i$, provided that
$D\in(0,D_*)$.
\end{lemma}
\proof  Let $E(y)$ be the error function
\begin{equation}\label{eqErrorFunction}
E(y)=\frac{2}{\sqrt{\pi}}\int_0^y \exp\left({-z^2}\right)\,dz.
\end{equation}
Then
\begin{equation}\label{eqPsiE}
\Psi(x,\tau)=\int_{x}^{\overline{x}} \psi(y,\tau)\,dy=E\left(
\frac{\bar x-x}{2\sqrt{\tau}}\right).
\end{equation}

For $i=1,2,\ldots$, we consider an increasing sequence of time
moments $t_i$, an (auxiliary) increasing sequence of times
$\theta_i$, and an increasing sequence $y_i$  such that
$t_{i-1}<\theta_i<t_i$ and
\begin{equation}\label{2ol}
E\left(\frac{y_{i-1}}{2\sqrt{t}}\right)<1/2,\quad t\ge\theta_i;
\end{equation}
\begin{equation}\label{1ol}
2E\left(\frac{y_i}{2\sqrt{t_i}}\right)-2E\left(\frac{y_{i-1}}{2\sqrt{\theta_i}}\right)-1>\frac{1}{\mu^2(t_i-\theta_i)},
\end{equation}
where $t_0=0$ and $y_0=0$. Such sequences can always be
constructed step by step.
For example, at the $i$th step, we can first fix a sufficiently
large $\theta_i$ to satisfy~\eqref{2ol}; then take a sufficiently
large $t_i>\theta_i$ to have
\begin{equation}\label{3olAux}
 2-2E\left(\frac{y_{i-1}}{2\sqrt{\theta_i}}\right)-1 >\frac{1}{\mu^2(t_i-\theta_i)};
\end{equation}
and then, using $E(\infty)=1$ and~\eqref{3olAux}, choose a
sufficiently large $y_i$ to ensure~\eqref{1ol}.
Since $E$ is an increasing function, relation \eqref{1ol} implies
\begin{equation}\label{3ol}
 2E\left(\frac{y_i}{2\sqrt{t}}\right)-2E\left(\frac{y_{i-1}}{2\sqrt{t}}\right)-1 >\frac{1}{\mu^2(t_i-\theta_i)},
\quad t\in [\theta_i,t_i],\
\end{equation}
for all $i=1,2,\dots$.

Now fix $N\in\bbN$. Using the relations
$\overline{x}-x_i=D^{1/2}y_i$, the equality~\eqref{eqPsiE}, and
Lemma~\ref{lphi-psi} for a specific $\theta=1$, we have
\begin{equation}\label{eqPhi-E}
\left|\Phi(x_i,Dt)-E\left(\frac{y_i}{2\sqrt{t}}\right)\right|\le
c_1 \exp\left({-\frac{(\ox-\ux)^2}{{5Dt}}}\right),\qquad
i=1,\ldots,N;\ t\in (0,t_N],
\end{equation}
where $c_1>0$ does not depend on $D$, provided that $D\in(0,D_*)$
and $D_*t_N\in(0,1]$.

Due to~\eqref{eqPhi-E}, by choosing $D_*$ small enough, we can
make the value
$\left|\Phi(x_i,Dt)-E\left(\frac{y_i}{2\sqrt{t}}\right)\right|$
arbitrarily small for all $i=1,\ldots,N$, $t\in (0,t_N]$, and
$D\in(0,D_*)$. Hence, using~\eqref{2ol} and~\eqref{3ol} and the
triangle inequality, we obtain inequalities~\eqref{3ol''}
and~\eqref{3ol'}.
\endproof

\begin{lemma}\label{lUclosePhiTheor2}
Let $D_*$ be the constant from
Lemma~$\ref{lPhiinequalitiesTheor2}$. Then, for any $D\in(0,D_*)$
and $\varepsilon_2>0$, there are
$\varepsilon_*=\varepsilon_*(\varepsilon_2,D)>0$ and
$v_*=v_*(\varepsilon_2)$ such that
$$
\left\|\dfrac{U(\cdot,t)}{\uU(t)} -
\Phi(\cdot,Dt)\right\|_{C[\underline{x},\overline{x}]}\le
 \varepsilon_2\quad \forall t\ge \theta_1,
$$
provided that   $u_0(x)$ satisfies~\eqref{data2} with
$\varepsilon\in(0,\varepsilon_*)$ and $v_0\in[0, v_*)$.
\end{lemma}
\proof In the proof of Lemma~\ref{lsmallutau0D}
(see~\eqref{equlowerupper} and~\eqref{equuoue}), we have shown
that
\begin{equation}\label{eqUclosePhiTheor2-1}
\int_x^\ox\uu(y,t)\,dy\le U(x,t)\le
\exp\left({\frac{v_0}{\mu}}\right)\int_x^\ox\uu(y,t)\,dy, \quad
x\in[\underline{x},\overline{x}],\ t\ge0,
\end{equation}
where $\uu$ solves
$$ \uu_t=D\uu_{xx}, \quad
\uu_x|_{x=\ox}=\uu_x|_{x=\ux}=0,\quad \uu|_{t=0}=u_0(x).
$$

On the other hand, Lemma~\ref{lphi-omega1} implies that
\begin{equation}\label{eqUclosePhiTheor2-2}
\begin{aligned}
& \left\|\uu(\cdot,t) -
\varphi(\cdot,Dt)\right\|_{C[\underline{x},\overline{x}]}\\
&\qquad \le C\left(\varepsilon^2 (D\theta_1)^{-3/2}+\varepsilon
(D\theta_1)^{-1/2}\right),\quad t\ge \theta_1,
\end{aligned}
\end{equation}
where $\theta_1$ is defined in Lemma~\ref{lPhiinequalitiesTheor2}
and $C>0$ does not depend on $t,\varepsilon,D,\theta_1$.

Now fix $D$ and an arbitrary $\varepsilon_1\in(0,1)$. Choose
$v_*=v_*(\varepsilon_1)$ such that
\begin{equation}\label{eqUclosePhiTheor2-3}
\exp\left({\frac{v_0}{\mu}}\right)\le 1+ \varepsilon_1,\quad
1-\varepsilon_1\le\dfrac{1}{\uU(t)}\le 1
\end{equation}
for all $v_0\in[0,v_*)$, where we use that $1=\uU(0)\le \uU(t)\le
\uU(0)+v(0)=1+v_0$. Next choose $\varepsilon_*$ such that
\begin{equation}\label{eqUclosePhiTheor2-4}
 \varepsilon^2 (D\theta_1)^{-3/2}+\varepsilon (D\theta_1)^{-1/2}\le
\varepsilon_1
\end{equation}
for all $\varepsilon\in(0,\varepsilon_*)$.

Combining~\eqref{eqUclosePhiTheor2-1}--\eqref{eqUclosePhiTheor2-4}
and using the fact that $|\Phi(x,t)|\le 1$ yields
$$
\left\|\dfrac{U(\cdot,t)}{\uU(t)} -
\Phi(\cdot,Dt)\right\|_{C[\underline{x},\overline{x}]}\le
 c_1\varepsilon_1\quad \forall t\ge \theta_1,
$$
where $c_1$ does not depend on the other quantities in the
inequality. Setting $\varepsilon_1=\varepsilon_2/c_1$, we complete
the proof.
\endproof

Now Lemma~\ref{lUinequalitiesTheor2} follows from
Lemmas~\ref{lPhiinequalitiesTheor2} and~\ref{lUclosePhiTheor2}.

\subsection{Proof of Theorem~\ref{tt1}}
Relations \eqref{5ol} and \eqref{4ol} are counterparts of
equations \eqref{2*} and \eqref{1*}, respectively. The proof can
be now completed by exactly the same argument as we used in
Sec.~\ref{subsecProofTheoremt1}. It suffices  to use estimates
\eqref{5ol} and \eqref{4ol}   in place of \eqref{2*} and
\eqref{1*}.

In particular, relation \eqref{5ol} ensures that every front
belonging to the interval $[x_{i-1},\overline{x})$ at the moment
$t=\theta_i$ is immortal and that if the interval $[x_i,x_{i-1})$
contains at least one front at the moment $t=\theta_i$, then the
most right of these fronts is immortal.

Relation \eqref{4ol} ensures that if the interval $[x_i,x_{i-1})$
does not contain fronts at the moment $t=\theta_i$, then at least
one immortal front arrives at this interval at a moment $t_*$
preceding $t_i$, because during the time interval $\theta_i\le t
\le \min\{t_*,t_i\}$, we have
$$
\dot w(t)\ge \mu^2 \uU(t) \left(2\frac{U(x_i,D
t)}{\uU(t)}-2\frac{U(x_{i-1},D
t)}{\uU(t)}-1\right)>\frac1{t_i-\theta_i}.
$$


\begin{remark}

We have also obtained asymptotic formulas (as the initial amount
of nutrient $v_0$ and the diffusion coefficient $D$ tend to zero,
while the initial density $u_0(x)$ tends to the delta function)
for positions of a few first steady immortal fronts and the time
moments at which these fronts appear.

Let us consider the formation of fronts on a fixed time interval
$[0,T)$, assuming that:
\begin{enumerate}
\item The initial density $u_0(x)$ satisfies~\eqref{data2} with
$\varepsilon=\varepsilon(D)>0$, where $\varepsilon(D)\to 0$ as
$D\to 0$.

\item The initial amount of food is small: $v_0=v_0(D)> 0$, where
$v_0(D)\to 0$  as $D\to 0$.

\item The initial value of the input is close to $\ox$:
$w_0=w_0(D)\le\ox$, where   $w_0(D)\to\ox$ as $D\to 0$.

\item The initial configuration is $r_0(x)\equiv 1$.
\end{enumerate}

\end{remark}

Under these assumptions, we denote by $t_1<t_2<\cdots$ the
consecutive time moments at which the moving fronts become steady
and immortal; their positions at these moments are denoted
$x_1>x_2>\cdots$ (where $x_i\in [\underline{x},\overline{x}]$ for
all $i$). We have obtained the formulas
$$
t_n=s_n+o(1);\qquad x_n= \ox-q_n ,\quad\text{where}\quad
q_n=2\big(y_n+o(1)\big)\big(Ds_n\big)^{1/2}
$$
for the first six fronts; here $o(1)$ stands for functions of $D$
that tend to $0$ as $D\to 0$ (the convergence of $o(1)$ to zero
may depend on $T$, but not on $n$). If, for example, $\ox=1/4$,
then the values of $s_n$, $y_n$ in these formulas (which, in fact,
do not depend on $\ux$) are equal to
$$
\begin{aligned}
&(s_1,y_1)\approx (2.2,\,0.48), & &(s_2,y_2)\approx
(9.1,\,0.83), & &(s_3,y_3)\approx (24.0,\,1.07),\\
&(s_4,y_4)\approx (53.3,\,1.26), & &(s_5,y_5)\approx
(108.4,\,1.42), & &(s_6,y_6)\approx (209.4,\,1.57),
\end{aligned}
$$
which implies that
$$
\begin{aligned}
&q_1\approx 1.4\,D^{1/2}, &  &q_2\approx 5.0\, D^{1/2}, &
&q_3\approx 10.5\,D^{1/2},\\
&q_4\approx 18.4\,D^{1/2}, & &q_5\approx 29.6\,D^{1/2}, &
&q_6\approx 45.4\,D^{1/2}.
\end{aligned}
$$
 The authors do not know asymptotics of further
fronts; this is a subject of future work.

\section*{Acknowledgments}
PG acknowledges the support of Collaborative Research
Center 910 (Germany). DR was supported by NSF grant DMS-1413223. The authors are grateful to Sergey
Tikhomirov who created a software for a number of numerical
experiments.






\begin{thebibliography}{10}


\expandafter\ifx\csname url\endcsname\relax
  \def\url#1{\texttt{#1}}\fi
\expandafter\ifx\csname
urlprefix\endcsname\relax\def\urlprefix{URL }\fi
\expandafter\ifx\csname href\endcsname\relax
  \def\href#1#2{#2} \def\path#1{#1}\fi

\bibitem{h2}
Applebe B, Flynn D, McNamara H, O'Kane P, Pimenov A, Pokrovskii A,
Rachinskii D, Zhezherun A (2009) Rate-independent hysteresis in
terrestrial hydrology. IEEE Control Systems Magazine 29 (1):
44-69.

\bibitem{h3}
Appelbe B, Rachinskii D, Zhezherun A (2008) Hopf bifurcation in a
van der Pol type oscillator with magnetic hysteresis. Physica B
403 (2-3): 301-304.

\bibitem{benzer}
Benzer S (1953) Induced synthesis of enzymes in bacteria analyzed
at the cellular level. Biochim et Biophys Acta 11: 383-395.

\bibitem{h5}
Brokate M, Pokrovskii AV, Rachinskii DI (2006) Asymptotic
stability of continual sets of periodic solutions to systems with
hysteresis. J Math Anal Appl 319: 94-109.

\bibitem{h1}
Brokate M, Pokrovskii AV, Rachinskii DI, Rasskazov O (2005)
Differential equations with hysteresis via a canonical example.
In: The Science of Hysteresis. Mayergoyz ID, Bertotti G, editor.,
Elsevier. Vol. I, Chapter II: 125-291.

\bibitem{cohn1}
Cohn M, Horibata K (1959) Inhibition by glucose of the induced
synthesis of the beta-galactoside-enzyme system of Escherichia
coli. Analysis of maintenance. J Bacteriol 78: 601-612.

\bibitem{cohn2}
Cohn M, Horibata K (1959) Analysis of the differentiation and of
the heterogeneity within a population of Eschericia coli
undergoing induced beta- galactosidase synthesis. J Bacteriol 78:
613-623.

\bibitem{cohn3}
Cohn M, Horibata K (1959) Physiology of the inhibition by glucose
of the induced synthesis of the beta-galactoside-enzyme system of
Escherichia coli. J Bacteriol 78: 624-635.

\bibitem{h4}
Cross R, McNamara H, Pokrovskii A, Rachinskii D (2008) A new
paradigm for modelling hysteresis in macroeconomic flows. Physica
B 403 (2-3): 231-236.

\bibitem{delbr}
    Delbr\"uck M (1949) Translation of discussion following a paper by Sonneborn TM and Beale GH. In: Unit\'e Biologiques Dou\'ees de Continuit\'e G\'en\'etique. Editions du Centre National de la Recherche Scientifique, Paris: 33-35.

\bibitem{h7}
Diamond P, Kuznetsov NA, Rachinskii DI (2001) On the Hopf
bifurcation in control systems with asymptotically homogeneous at
infinity bounded nonlinearities. J Differential~Equations 175:
1-26.

\bibitem{h6}
Diamond P, Rachinskii DI, Yumagulov MG (2000) Stability of large
cycles in a nonsmooth problem with Hopf bifurcation at infinity.
Nonlinear Anal 42(6): 1017-1031.

\bibitem{c10}
Dubnau D, Losick R (2006) Bistability in bacteria. Mol Microbiol
61: 564-572.

\bibitem{1}
Friedman G, Gurevich P, McCarthy S, Rachinskii D (2014) Switching
behaviour of two-phenotype bacteria in varying environment. J Phys
Conf Ser, accepted.

\bibitem{FD}
Friedman G, Rachinskii D (2014) Hysteresis can grant fitness in stochastically varying environment.
PLoS ONE 9(7): e103241.

\bibitem{15}   Gardner TS, Cantor CR, Collins JJ (2000) Construction of a genetic toggle switch in Escherichia coli. Nature 403: 339-342.

\bibitem{grazi}
Graziani S, Silar P, Daboussi M-J (2004) Bistability and
hysteresis of the `Secteur' differentiation are controlled by a
two-gene locus in Nectria haematococca. BMC Biology 2 (18): 1-19.

\bibitem{2}
Gurevich P, Rachinskii D (2013) Well-posedness of parabolic
equations containing hysteresis with diffusive thresholds.
Proc Steklov Inst Math 283: 87-109.

\bibitem{GurevichTikhomirovEquadiff13}
Gurevich PL, Tikhomirov SB (2014) Systems of reaction-diffusion
equations with spatially distributed hysteresis. Mathematica
Bohemica (Proc Equadiff 2013), accepted.

\bibitem{ark}
Ham TS, Lee SK, Keasling JD, Arkin AP (2008) Design and
construction of a double inversion recombination switch for
heritable sequential genetic memory. PLoS ONE 3 (7): e2815.

\bibitem{Jaeger1}
Hoppensteadt FC, J\"ager W (1980) Pattern formation by bacteria.
Lecture Notes in Biomathematics 38: 68-81.

\bibitem{Kopfova}
Kopfova J (2006) Hysteresis in biological models. J Phys Conf Ser
55:  130-134.

\bibitem{h10}
Krasnosel'skii AM, Rachinskii DI (2001) On continua of cycles in
systems with hysteresis. Doklady Math 63 (3): 339-344.

\bibitem{h9}
Krasnosel'skii AM, Rachinskii DI (2002) On a bifurcation governed
by hysteresis nonlinearity. NoDEA Nonlinear Differential Equations
Appl 9: 93-115.

\bibitem{h8}
Krejci P, O'Kane P, Pokrovskii A, Rachinskii D (2011) Stability
results for a soil model with singular hysteretic hydrology. J
Phys Conf Ser 268: 012016.

\bibitem{bac2}
Kussell E, Lieber S, Phenotypic diversity, population growth, and
information in fluctuating environments. Science 309: 2075-2078.

\bibitem{lai}
Lai K, Robertson MJ, Schaffer DV (2004) The sonic hedgehog
signaling system as a bistable genetic switch. Biophys J 86:
2748-2757.

\bibitem{SH}
(2005) The Science of Hysteresis. Mayergoyz ID, Bertotti G,
editors. Elsevier. 751 p.

\bibitem{maamar}
Maamar H, Raj A, Dubnau D (2007) Noise in gene expression
determines cell fate in Bacillus subtilis. Science 317: 526-529.

\bibitem{monod}
Monod J (1966) From enzymatic adaptation to allosteric
transitions. Science 154: 475-483.

\bibitem{f0} Novick A, Weiner M (1957)
Enzyme induction as an all-or-none phenomenon. Proc Natl Acad Sci
USA 43 (7): 553-566.

\bibitem{oud}
Ozbudak EM, Thattai M, Lim HN, Shraiman BI, van Oudenaarden A
(2004) Multistability in the lactose utilization network of
Escherichia coli. Nature 427: 737-740.

\bibitem{10}
    Pomerening JR, Sontag ED, Ferrell JE Jr (2003) Building a cell cycle oscillator: hysteresis and bistability in the activation of Cdc2. Nature Cell Biol 5: 346-351.

\bibitem{smits}
Smits WK, Kuipers OP, Veening JW (2006) Phenotypic variation in
bacteria: the role of feedback regulation. Nat Rev Microbiol 4:
259-271.

\bibitem{Spiegelman}
Spiegelman S, DeLorenzo WF (1952) Substrate stabilization of
enzymeforming capacity during the segregation of a heterozygote.
Proc Natl Acad Sci U S A 38 (7): 583-592.

\bibitem{bac1}
Thattai M, van Oudenaarden A (2004) Stochastic gene expression in
fluctuating environments. Genetics 167: 523-530.

\bibitem{f3} Thattai M, Shraiman B (2003)
Metabolic switching in the sugar phosphotransferase system of
Escherichia coli. Biophys J 85: 744-754.

\bibitem{VisintinSpatHyst}
Visintin A (1986) Evolution problems with hysteresis in the source
term. SIAM J Math Anal 17: 1113-1138.

\bibitem{wanga}
Wanga L, Walkera BL, Iannacconeb S, Bhatta D, Kennedya PJ, Tse WT
(2009) Bistable switches control memory and plasticity in cellular
differentiation. PNAS U S A 106 (16): 6638-6643.

\bibitem{winge}
Winge \"O, Roberts C (1948) Inheritance of enzymatic characters in
yeast and the phenomenon of long term adaptation. Comp rend trav
Lab, Carlsberg. Ser physiol 24: 263-315.

\bibitem{c12}
Wolf DM, Arkin AP (2003) Motifs, modules and games in bacteria.
Curr Opin Microbiol: 125-134.

\bibitem{bac4}
Wolf DM, Fontaine-Bodin L, Bischofs I, Price G, Keasling J, Arkin
AP (2008) Memory in microbes: quantifying history-dependent
behaviour in a bacterium. PLoS ONE 3 (2): e1700.

\bibitem{h2}
Wolf DM, Vazirani VJ, Arkin AP (2005) Diversity in times of
adversity: probabilistic strategies in microbial survival games. J
Theor Biol 234 (2): 227-253.

\end{thebibliography}







\end{document}